\documentclass[11pt,a4paper]{article}

\usepackage{inputenc}
\usepackage{amsmath}
\usepackage{bm}
\usepackage{bbold}
\usepackage{amsthm}
\usepackage{enumerate}

\usepackage{graphicx}
\usepackage{fontawesome}

\usepackage{hyperref}

\title{Using tropical optimization to solve constrained minimax single-facility location problems with rectilinear distance}

\author{N. Krivulin\thanks{Faculty of Mathematics and Mechanics, Saint Petersburg State University, 28 Universitetsky Ave., St.~Petersburg, 198504, Russia, 
nkk@math.spbu.ru.}
}

\date{}

\newtheorem{theorem}{Theorem}

\newtheorem{corollary}[theorem]{Corollary}

\theoremstyle{definition}

\begin{document}

\maketitle

\begin{abstract}
The aim of this paper is twofold: first, to extend the area of applications of tropical optimization by solving new constrained location problems, and second, to offer new closed-form solutions to general problems that are of interest to location analysis. We consider a constrained minimax single-facility location problem with addends on the plane with rectilinear distance. The solution commences with the representation of the problem in a standard form, and then in terms of tropical mathematics, as a constrained optimization problem. We use a transformation technique, which can act as a template to handle optimization problems in other application areas, and hence is of independent interest. To solve the constrained optimization problem, we apply methods and results of tropical optimization, which provide direct, explicit solutions. The results obtained serve to derive new solutions of the location problem, and of its special cases with reduced sets of constraints, in a closed form, ready for practical implementation and immediate computation. As illustrations, numerical solutions of example problems and their graphical representation are given. We conclude with an application of the results to optimal location of the central monitoring facility in an indoor video surveillance system in a multi-floor building environment.
\\

\textbf{Key-Words:} minimax location problem, rectilinear distance, idempotent semifield, tropical optimization, constrained optimization, explicit solution.
\\

\textbf{MSC (2010):} 90B85, 15A80, 65K05, 90C48
\end{abstract}

\section{Introduction}

Tropical (idempotent) mathematics, originated in the middle of the last century as the theory and applications of semirings with idempotent addition, finds use in a variety of fields, from operations research to algebraic geometry. The significant advances, achieved in the area of tropical mathematics in the last decades, are reported in several monographs, including recent ones by \cite{Golan2003Semirings,Heidergott2006Maxplus,Mceneaney2006Maxplus,Itenberg2007Tropical,Gondran2008Graphs,Maclagan2015Introduction}, and in a wide range of research papers.  

Optimization problems that are formulated and solved in terms of tropical mathematics are a matter of concern for tropical optimization, which presents an important research domain, with the focus on new solutions to old and fresh problems in operations research and management science. Applications of tropical optimization include real-world problems in project scheduling, location analysis, transportation networks, discrete event dynamic systems, decision making, and in other fields. 

Location problems constitute one of the classical areas in optimization, which dates back to the XVII century. A variety of approaches and techniques exists to solve location problems in different settings, including methods of mathematical programming, and of discrete, combinatorial and graph optimization (see, e.g. \cite{Sule2001Logistics,Klamroth2002Singlefacility,Farahani2009Facility,Eiselt2011Foundations,Laporte2015Location} for the current state of the art in the area). 

There are certain location problems that have solutions obtained in the framework of tropical optimization. Specifically, a solution in terms of tropical mathematics is proposed by \cite{Cuninghamegreen1991Minimax,Cuninghamegreen1994Minimax} to one-dimensional minimax location problems defined on graphs. Furthermore, several constrained minimax location problems are examined by \cite{Zimmermann1991Minmax,Zimmermann1992Optimization,Hudec1993Aservice,Hudec1999Biobjective,Tharwat2010Oneclass} in the context of the theory of max-separable functions, which is closely related to the tropical mathematics approach. Finally, methods of tropical optimization are applied to solve unconstrained and constrained minimax single-facility location problems with Chebyshev and rectilinear distances \cite{Krivulin2011Algebraic,Krivulin2011Anextremal,Krivulin2012Anew,Krivulin2013Direct,Krivulin2014Complete,Krivulin2015Onanalgebraic}.

The aim of this paper is twofold: first, to develop new applications of tropical optimization by solving new location problems, and second, to offer new closed-form solutions to rather general problems that are of interest to location analysis. We consider a constrained minimax single-facility location problem with addends on the plane with rectilinear distance, which can be referred to as a constrained Rawls location problem or a constrained messenger boy problem. The solution commences with the representation of the problem, first formulated in a standard form, in terms of tropical mathematics as a constrained optimization problem. We use a transformation technique, which can act as a template to handle optimization problems in other application areas, and hence is of independent interest. To solve the constrained optimization problem, we apply methods and results of tropical optimization \cite{Krivulin2014Complete,Krivulin2015Extremal,Krivulin2015Amultidimensional,Krivulin2017Direct}, which provide direct, explicit solutions of the problem and of its special cases with reduced sets of constraints. We further develop the methods to extend the solution of the unconstrained problem provided by \cite{Krivulin2011Anextremal,Krivulin2015Onanalgebraic} to the constrained problems of interest. The results are obtained in a closed form, ready for immediate computation, and can serve to complement and supplement known solutions of the location problems under examination. To illustrate application of the results, we describe a direct solution to the problem of optimal location of the central monitoring facility in an indoor closed-circuit television (CCTV) video surveillance system in a multi-floor building environment.

The rest of the paper is organized as follows. We begin with Section~\ref{S-CMRSFLP}, where the location problem of interest is formulated in a standard way. Section~\ref{S-PADN} includes an overview of the definitions and notation of idempotent algebra to be used in the subsequent sections. In Section~\ref{S-TOP}, we consider several tropical optimization problems and describe their solutions. Section~\ref{S-TSLP} offers the main result of the paper. First, we represent the location problem under consideration as a constrained tropical optimization problem, and then solve this optimization problem using the results of the previous section. As a consequence, solutions are obtained to some special cases of the problem with reduced sets of constraints. We use the results given in terms of tropical mathematics to derive solutions of the location problems in the standard form. In Section~\ref{S-NEGI}, we present numerical examples and offer graphical illustrations. We conclude with an application to an optimal location problem, arising in the deployment of an indoor CCTV video surveillance system, in Section~\ref{S-ACCTVMFL}, and make some final observations and comments in Section~\ref{S-C}.

\section{Constrained minimax rectilinear single-facility location problem}
\label{S-CMRSFLP}

We start with a brief outline of the optimization problem, which is drawn from location analysis to motivate the present study. A comprehensive overview of various location problems, their solutions and application examples is provided by a series of surveys published at different times, including \cite{Francis1983Locational,Brandeau1989Anoverview,ReVelle2005Location,Brimberg2009Optimizing,Chhajed2013Facility}. Further details can be found in monographs and collections of studies, such as recent books by \cite{Sule2001Logistics,Klamroth2002Singlefacility,Farahani2009Facility,Eiselt2011Foundations,Laporte2015Location}.

We consider a quite general problem to locate a new point (a facility center) on the plane to minimize the maximum of rectilinear distances to given points (demand centers), each of which can be modified by adding a constant called the addend. The optimal location is subject to constraints that impose upper bounds on distances from each given point to the new point, and define a strip-shaped feasible location region.  

The rectilinear metric (also known as the rectangular, Manhattan, right-angle, city-block, taxicab or $L_{1}$ metric) arises in location analysis in various applied contexts. Examples include locating a public or commercial facility in an urban area with a grid of rectangular streets, an industrial facility within a plant or warehouse with a system of perpendicular transport aisles, and an electronic component on an integrated circuit with orthogonal mesh of wires. The addends can represent an additional cost or distance required to reach each demand point, such as vertical distance when the rectilinear metric is defined on the horizontal plane.

Constrained minimax location problems appear in a range of application areas from urban planning to industrial and electrical engineering. A typical example is the optimal location of emergency service facilities (hospitals, police and fire stations, emergency shelters) in urban design, under constraints on the travel distances prescribed by emergency service standards and rules set by federal, state or municipal agencies. Since the minimax objectives in locating public facilities can well be interpreted in the framework of the theory of justice of John Rawls, these problems are frequently referred to as the Rawls location problems \cite{Hansen1981Constrained,Hansen1981Outcomes}. In addition, minimax single-facility location problems with and without addends are sometimes called the messenger boy problems and the delivery boy problems, respectively \cite{Elzinga1972Geometrical}.

We now represent the location problem under study in a formal way. First note that the rectilinear distance between vectors $\bm{a}=(a_{1},a_{2})^{T}$ and $\bm{b}=(b_{1},b_{2})^{T}$ in the real plane $\mathbb{R}^{2}$ is calculated as
\begin{equation*}
\rho(\bm{a},\bm{b})
=
|a_{1}-b_{1}|+|a_{2}-b_{2}|.
\end{equation*}

Suppose there is a set of $m\geq1$ given points, denoted by $\bm{r}_{j}=(r_{1j},r_{2j})^{T}\in\mathbb{R}^{2}$ for all $j=1,\ldots,m$. Let $w_{j}\in\mathbb{R}$ be the addend, associated with point $j$, and $d_{j}\in\mathbb{R}$, where $d_{j}\geq0$, be the upper bound on the distance to point $j$. Let $s,t\in\mathbb{R}$, where $s\leq t$, be the left and right boundary of the vertical strip, representing the feasible location area.

Then, the problem of interest, which can be referred to as the constrained minimax rectilinear single-facility location problem with addends, is formulated to find points $\bm{x}=(x_{1},x_{2})^{T}\in\mathbb{R}^{2}$ that  
\begin{equation*}
\begin{aligned}
&
\text{minimize}
&&
\max_{1\leq j\leq m}(\rho(\bm{x},\bm{r}_{j})+w_{j}),
\\
&
\text{subject to}
&&
\rho(\bm{x},\bm{r}_{j})
\leq
d_{j},
\quad
j=1,\ldots,m;
\\
&&&
s
\leq
x_{1}
\leq
t.
\end{aligned}
\end{equation*}

After rewriting the rectilinear distance $\rho$ in coordinate form, the problem becomes
\begin{equation*}
\begin{aligned}
&
\text{minimize}
&&
\max_{1\leq j\leq m}(|x_{1}-r_{1j}|+|x_{2}-r_{2j}|+w_{j}),
\\
&
\text{subject to}
&&
|x_{1}-r_{1j}|+|x_{2}-r_{2j}|
\leq
d_{j},
\quad
j=1,\ldots,m;
\\
&&&
s
\leq
x_{1}
\leq
t.
\end{aligned}
\end{equation*}

Consider the constraints in the problem. For each $j=1,\ldots,m$, the inequality $|x_{1}-r_{1j}|+|x_{2}-r_{2j}|\leq d_{j}$ defines on the plane a square rotated by $45^{\circ}$ around its center at the point $\bm{r}_{j}=(r_{1j},r_{2j})^{T}$, which is often called the diamond. The common area of all inequalities, if it exists, takes the form of a rectangle tilted $45^{\circ}$ to the axes. The feasible location region is the intersection of this rectangle with the strip area given by the inequality $s\leq x_{1}\leq t$, provided that the intersection is not empty.

Both special cases and extensions of the rectilinear single-facility location problem are thoroughly examined in the literature. For some unconstrained versions of the problem, direct solutions are obtained in a closed form, whereas, in other cases, the problems have solutions given by iterative computational algorithms, which find a solution, if it exists, or indicate that there are no solutions. Specifically, the unconstrained problem is considered without addends in \cite{Francis1972Ageometrical,Elzinga1972Geometrical}, and with addends in \cite{Elzinga1972Geometrical}, where closed-form solutions are derived based on geometric arguments.

An algorithmic solution is proposed by \cite{Dearing1972Onsome} for a weighted extension of the problem, in which the distances appear in the objective function with non-negative weights, and for a weighted multi-facility problem. A weighted multi-facility problem without constraints is solved by means of linear programming computational schemes in \cite{Wesolowsky1972Rectangular,Elzinga1973Anote}, whereas the constrained problem is by a network flow algorithm in \cite{Dearing1974Anetwork}. An interactive computer graphical technique is developed in \cite{Brady1980Interactive} to solve single-facility location problems with non-convex feasible regions.

An algebraic approach, which uses results of tropical optimization, is applied in \cite{Krivulin2011Anextremal,Krivulin2015Onanalgebraic} to the problem under consideration when all constraints are removed. The approach offers a direct, explicit solution based on a straightforward algebraic technique, rather than on geometric considerations in the classical works \cite{Francis1972Ageometrical,Elzinga1972Geometrical}.

Below, we further develop the algebraic approach to extend methods of tropical optimization to the constrained location problem with rectilinear distance. Based on this approach, we derive closed-form solutions for the location problem of interest, as well as for its special cases with reduced sets of constraints. To the best of our knowledge, no direct solutions to the location problem have been previously reported.

To handle the problem, we first transform it into a tropical optimization problem examined in \cite{Krivulin2014Complete} in the context of constrained location with Chebyshev distances ($L_{\infty}$ metric) to exploit the complete solution derived therein. Note that, from the geometrical point of view, the location problems on the plane with rectilinear and Chebyshev distances are known to convert into each other by rotation of the coordinate axes (see, e.g., \cite{Farahani2009Facility}), which can serve as additional intuition and evidence in support of the algebraic technique proposed for the transformation. Next, we further develop and refine obtained results to provide complete solutions of the location problems of interest, which are then represented both in terms of tropical mathematics and in the conventional form. 

To conclude this section, we observe that it is not difficult to represent the problem under study as a linear program, and then to solve it using methods and computational techniques of linear programming. However, these methods normally offer algorithmic solutions, and do not guarantee a direct solution in a closed form.

\section{Preliminary algebraic definitions and notation}
\label{S-PADN}

We now present a concise overview of the definitions and notation of idempotent algebra from \cite{Krivulin2014Complete,Krivulin2015Extremal,Krivulin2015Amultidimensional,Krivulin2017Direct}, which provide the basis for the description of the tropical optimization problems and their solutions in the next section, and for the use of these solutions to attack location problems in the subsequent sections. Further details on tropical mathematics are available in various introductory and advanced texts, including \cite{Golan2003Semirings,Heidergott2006Maxplus,Mceneaney2006Maxplus,Itenberg2007Tropical,Gondran2008Graphs,Maclagan2015Introduction} to name only a few recent publications.

\subsection{Idempotent semifield}

Let $\mathbb{X}$ be a non-empty set that is closed under two associative and commutative operations, addition $\oplus$ and multiplication $\otimes$, and equipped with two distinct elements, zero $\mathbb{0}$ and one $\mathbb{1}$, which are neutral with respect to addition and multiplication. Addition is idempotent, which indicates that the equality $x\oplus x=x$ holds for each $x\in\mathbb{X}$. Multiplication is distributive over addition, and invertible, which provides each $x\ne\mathbb{0}$ with an inverse $x^{-1}$ such that $x\otimes x^{-1}=\mathbb{1}$. Under these assumptions, the algebraic system $(\mathbb{X},\mathbb{0},\mathbb{1},\oplus,\otimes)$ is frequently called the idempotent semifield.

Idempotent addition provides a partial order on $\mathbb{X}$ to define $x\leq y$ if and only if $x\oplus y=y$. It follows from the definition of the order relation that the operations in the semifield have the following properties. First, the inequality $x\oplus y\leq z$ is equivalent to the two inequalities $x\leq z$ and $y\leq z$ for any $x,y,z\in\mathbb{X}$. Furthermore, both addition and multiplication are isotone, which implies that the inequality $x\leq y$ results in the inequalities $x\oplus z\leq y\oplus z$ and $x\otimes z\leq y\otimes z$. Finally, inversion is antitone, which means that $x\leq y$ yields $x^{-1}\geq y^{-1}$, provided that $x\ne\mathbb{0}$ and $y\ne\mathbb{0}$. In addition, the set $\mathbb{X}$ is assumed totally ordered by a linear order relation that is consistent with the partial order associated with addition.

As usual, the multiplication sign $\otimes$ is omitted below to safe writing. The integer power notation serves to signify iterated products, defined as $x^{0}=\mathbb{1}$, $x^{p}=xx^{p-1}$, $x^{-p}=(x^{-1})^{p}$ and $\mathbb{0}^{p}=\mathbb{0}$ for all $x\ne\mathbb{0}$ and integer $p>0$. The power notation is assumed extendable to allow rational and real exponents. 

A representative example of the idempotent semifield under consideration is the real semifield $\mathbb{R}_{\max,+}=(\mathbb{R}\cup\{-\infty\},-\infty,0,\max,+)$, where addition is defined as the operation of taking the maximum and multiplication is as the arithmetic addition, whereas the zero is given by $-\infty$ and the one is by $0$. For each $x\in\mathbb{R}$, there exists the inverse $x^{-1}$, which coincides with $-x$ in conventional algebra. The power $x^{y}$ acts as the arithmetic product $xy$ defined for any $x,y\in\mathbb{R}$. The partial order, which is given by addition, conforms to the standard linear order defined on $\mathbb{R}$. Finally, the obvious equality $\min(x,y)=-\max(-x,-y)$ yields $\min(x,y)=(x^{-1}\oplus y^{-1})^{-1}$ for all $x,y\in\mathbb{R}$.

As another example, consider the semifield $\mathbb{R}_{\min,\times}=(\mathbb{R}_{+}\cup\{+\infty\},+\infty,1,\min,\times)$, where $\mathbb{R}_{+}$ is the set of positive reals. This semifield has $\oplus=\min$, $\otimes=\times$, $\mathbb{0}=+\infty$ and $\mathbb{1}=1$. Both inversion and exponentiation have standard interpretation. The partial order induced by idempotent addition is opposite to the natural order on $\mathbb{R}$.

\subsection{Vector and matrix algebra}

The scalar addition $\oplus$ and multiplication $\otimes$ defined on $\mathbb{X}$ are routinely extended to vector and matrix operations. The set of matrices with $m$ rows and $n$ columns over $\mathbb{X}$ is denoted $\mathbb{X}^{m\times n}$. Addition and multiplication of matrices, and multiplication by scalars follow the standard rules. For any matrices $\bm{A}=(a_{ij})\in\mathbb{X}^{m\times n}$, $\bm{B}=(b_{ij})\in\mathbb{X}^{m\times n}$ and $\bm{C}=(c_{ij})\in\mathbb{X}^{n\times l}$, and a scalar $x\in\mathbb{X}$, the matrix operations are defined by the entry-wise formulae
$$
\{\bm{A}\oplus\bm{B}\}_{ij}
=
a_{ij}\oplus b_{ij},
\qquad
\{\bm{A}\bm{C}\}_{ij}
=
\bigoplus_{k=1}^{n}a_{ik}c_{kj},
\qquad
\{x\bm{A}\}_{ij}
=
xa_{ij}.
$$

The partial order relation and its associated properties of the operations in $\mathbb{X}$ extend entry-wise to the matrix operations. 

Consider square matrices of order $n$, which form the set $\mathbb{X}^{n\times n}$. A matrix whose diagonal entries are all equal to $\mathbb{1}$, and off-diagonal entries are to $\mathbb{0}$ is the identity matrix denoted by $\bm{I}$. The power notation is defined as $\bm{A}^{0}=\bm{I}$, $\bm{A}^{p}=\bm{A}\bm{A}^{p-1}$ for any square matrix $\bm{A}$ and integer $p>0$ to indicate repeated multiplication.

The trace of a matrix $\bm{A}=(a_{ij})\in\mathbb{X}^{n\times n}$ is calculated as $\mathop\mathrm{tr}\bm{A}=a_{11}\oplus\cdots\oplus a_{nn}$. The traces of matrix powers from $1$ to $n$ combine together to define the scalar
$$
\mathop\mathrm{Tr}(\bm{A})
=
\mathop\mathrm{tr}\bm{A}\oplus\cdots\oplus\mathop\mathrm{tr}\bm{A}^{n}.
$$

Provided that $\mathop\mathrm{Tr}(\bm{A})\leq\mathbb{1}$, the asterate operator (also known as the Kleene star) maps the matrix $\bm{A}$ to the matrix
$$
\bm{A}^{\ast}
=
\bm{I}\oplus\bm{A}\oplus\cdots\oplus\bm{A}^{n-1}.
$$

Any matrix that consists of one row (column) specifies a row (column) vector. All vectors below are assumed column vectors unless otherwise specified. The set of all column vectors with $n$ elements over $\mathbb{X}$ is denoted $\mathbb{X}^{n}$. 

A vector with all elements equal to $\mathbb{0}$ is the zero vector. Any vector without zero elements is called regular. For any vectors $\bm{a}=(a_{i})$ and $\bm{b}=(b_{i})$ in $\mathbb{X}^{n}$, and a scalar $x\in\mathbb{X}$, the vector addition and scalar multiplication are given by
$$
\{\bm{a}\oplus\bm{b}\}_{i}
=
a_{i}\oplus b_{i},
\qquad
\{x\bm{a}\}_{i}
=
xa_{i}.
$$

For any nonzero vector $\bm{a}=(a_{i})\in\mathbb{X}^{n}$, the multiplicative conjugate transpose is a row vector $\bm{a}^{-}=(a_{i}^{-})$ that has elements $a_{i}^{-}=a_{i}^{-1}$ if $a_{i}\ne\mathbb{0}$, and $a_{i}^{-}=\mathbb{0}$ otherwise.

The conjugate transposition is antitone in the sense that, if regular vectors $\bm{a}$ and $\bm{b}$ satisfy the element-wise inequality $\bm{a}\leq\bm{b}$, then $\bm{a}^{-}\geq\bm{b}^{-}$. In addition, the transposition has the following properties. For any nonzero vector $\bm{a}$, the equality $\bm{a}^{-}\bm{a}=\mathbb{1}$ holds. If the vector $\bm{a}$ is regular, then the entry-wise inequality $\bm{a}\bm{a}^{-}\geq\bm{I}$ is also valid.

\section{Tropical optimization problems}
\label{S-TOP}

In this section, we use idempotent algebra to formulate optimization problems and to describe their solutions. The problems find applications in various fields, including location analysis. Specifically, such problems occur in solving unconstrained and constrained single-facility location problems in the multidimensional space with Chebyshev distance \cite{Krivulin2012Anew,Krivulin2013Direct,Krivulin2014Complete}.

In the succeeding sections, we extend the solutions presented here to constrained single-facility location problems defined on the plane with rectilinear metric.  

We start with a general constrained optimization problem formulated in terms of an arbitrary idempotent semifield. Suppose that, given vectors $\bm{p},\bm{q},\bm{g},\bm{h}\in\mathbb{X}^{n}$, and a matrix $\bm{B}\in\mathbb{X}^{n\times n}$, the problem is to find all regular vectors $\bm{x}\in\mathbb{X}^{n}$ that
\begin{equation}
\begin{aligned}
&
\text{minimize}
&&
\bm{x}^{-}\bm{p}\oplus\bm{q}^{-}\bm{x},
\\
&
\text{subject to}
&&
\bm{B}\bm{x}
\leq
\bm{x},
\\
&
&&
\bm{g}
\leq
\bm{x}
\leq
\bm{h}.
\end{aligned}
\label{P-xpqx-Bxx-gxh}
\end{equation}

The solution of the problem, given in \cite{Krivulin2014Complete}, involves the introduction of a parameter to represent the minimum value of the objective function. Then, the problem reduces to solving a parametrized system of linear inequalities. The existence conditions of regular solutions for the system serve to evaluate the parameter, whereas the solution of the system is taken as a complete solution to the initial optimization problem. The results obtained take the form of the following statement.
\begin{theorem}
\label{T-xpqx-Bxx-gxh}
Let $\bm{B}$ be a matrix with $\mathop\mathrm{Tr}(\bm{B})\leq\mathbb{1}$, $\bm{p}$ be a nonzero vector, $\bm{q}$ and $\bm{h}$ be regular vectors, and $\bm{g}$ be a vector such that $\bm{h}^{-}\bm{B}^{\ast}\bm{g}\leq\mathbb{1}$. 

Then, the minimum value in problem \eqref{P-xpqx-Bxx-gxh} is equal to
\begin{equation}
\theta
=
(\bm{q}^{-}\bm{B}^{\ast}\bm{p})^{1/2}
\oplus
\bm{h}^{-}\bm{B}^{\ast}\bm{p}\oplus\bm{q}^{-}\bm{B}^{\ast}\bm{g},
\label{E-theta-qBp-hBpqBg}
\end{equation}
and all regular solutions of the problem are given by
$$
\bm{x}
=
\bm{B}^{\ast}\bm{u},
$$
where $\bm{u}$ is any regular vector such that
\begin{equation}
\bm{g}\oplus\theta^{-1}\bm{p}
\leq
\bm{u}
\leq
((\bm{h}^{-}\oplus\theta^{-1}\bm{q}^{-})\bm{B}^{\ast})^{-}.
\label{E-gthetap-u-hthetaqB}
\end{equation}
\end{theorem}

The conditions of the theorem have the following meaning. The requirement of a regular vector $\bm{h}$ is a necessary condition that allows regular solutions of the inequality $\bm{g}\leq\bm{x}\leq\bm{h}$. The condition $\mathop\mathrm{Tr}(\bm{B})\leq\mathbb{1}$ is necessary and sufficient to have a set of regular solutions to the inequality $\bm{B}\bm{x}\leq\bm{x}$, whereas $\bm{h}^{-}\bm{B}^{\ast}\bm{g}\leq\mathbb{1}$ is for a non-empty intersection of this solution set with the set defined by the constraints $\bm{g}\leq\bm{x}\leq\bm{h}$.

The assumptions of a non-zero vector $\bm{p}$ and a regular vector $\bm{q}$ are sufficient to keep the minimum value $\theta>\mathbb{0}$, which allows the inverse $\theta^{-1}$ to exist. These two assumptions can be replaced by a list of weaker conditions, which is, however, insufficient for application to the location problems under consideration. 

Consider two consequences of the theorem, which solve problem \eqref{P-xpqx-Bxx-gxh} when one of the inequality constraints is eliminated. First, we exclude the double inequality to write the problem
\begin{equation}
\begin{aligned}
&
\text{minimize}
&&
\bm{x}^{-}\bm{p}\oplus\bm{q}^{-}\bm{x},
\\
&
\text{subject to}
&&
\bm{B}\bm{x}
\leq
\bm{x}.
\end{aligned}
\label{P-xpqx-Bxx}
\end{equation}

The general solution, which is offered by Theorem~\ref{T-xpqx-Bxx-gxh}, takes the form of the next result (see, also \cite{Krivulin2012Anew}).
\begin{corollary}
\label{C-xpqx-Bxx}
Let $\bm{B}$ be a matrix with $\mathop\mathrm{Tr}(\bm{B})\leq\mathbb{1}$, $\bm{p}$ be a non-zero vector, and $\bm{q}$ be a regular vector.
Then, the minimum value in problem \eqref{P-xpqx-Bxx} is equal to
\begin{equation*}
\theta
=
(\bm{q}^{-}\bm{B}^{\ast}\bm{p})^{1/2},
\end{equation*}
and all regular solutions are given by
\begin{equation*}
\bm{x}
=
\bm{B}^{\ast}\bm{u},
\qquad
\theta^{-1}\bm{p}
\leq
\bm{u}
\leq
\theta(\bm{q}^{-}\bm{B}^{\ast})^{-}.
\end{equation*}
\end{corollary}

Furthermore, we consider another special case of the constrained problem at \eqref{P-xpqx-Bxx-gxh}, formulated to
\begin{equation}
\begin{aligned}
&
\text{minimize}
&&
\bm{x}^{-}\bm{p}\oplus\bm{q}^{-}\bm{x},
\\
&
\text{subject to}
&&
\bm{g}
\leq
\bm{x}
\leq
\bm{h}.
\end{aligned}
\label{P-xpqx-gxh}
\end{equation}

A solution goes as follows (see, also \cite{Krivulin2013Direct}).
\begin{corollary}
\label{C-xpqx-gxh}

Let $\bm{p}$ be a non-zero vector, $\bm{q}$ and $\bm{h}$ be regular vectors, and $\bm{g}$ be a vector such that $\bm{g}\leq\bm{h}$. Then, the minimum value in problem \eqref{P-xpqx-gxh} is equal to
\begin{equation*}
\theta
=
(\bm{q}^{-}\bm{p})^{1/2}\oplus\bm{h}^{-}\bm{p}\oplus\bm{q}^{-}\bm{g},
\end{equation*}
and all regular solutions of the problem are given by the condition
$$
\bm{g}\oplus\theta^{-1}\bm{p}
\leq
\bm{x}
\leq
(\bm{h}^{-}\oplus\theta^{-1}\bm{q}^{-})^{-}.
$$
\end{corollary}

Finally, we present a solution to the unconstrained problem \cite{Krivulin2012Anew}
\begin{equation}
\begin{aligned}
&
\text{minimize}
&&
\bm{x}^{-}\bm{p}\oplus\bm{q}^{-}\bm{x}.
\end{aligned}
\label{P-xpqx}
\end{equation}

\begin{corollary}
\label{C-xpqx}
Let $\bm{p}$ be a non-zero vector, and $\bm{q}$ be a regular vector. Then, the minimum value in problem \eqref{P-xpqx} is equal to
\begin{equation*}
\theta
=
(\bm{q}^{-}\bm{p})^{1/2},
\end{equation*}
and all regular solutions are given by
\begin{equation*}
\theta^{-1}\bm{p}
\leq
\bm{x}
\leq
\theta\bm{q}.
\end{equation*}
\end{corollary}

Below we show how to apply the results of tropical optimization in this section to solve the constrained rectilinear single-facility location problem under study.

\section{Transformation and solution of location problem}
\label{S-TSLP}

We are now in a position to turn back to the location problem formulated above. The solution begins with the representation of the problem in terms of tropical mathematics. We apply a useful transformation technique, which can serve as a template to handle other optimization problems, and hence is of independent interest. Furthermore, we derive, in the framework of tropical optimization, a direct solution to the problem in the general setting, and then present solutions to special cases of the problem, where some or all of constraints are removed. The section concludes with direct representation of the solutions in terms of the conventional algebra.

\subsection{Representation in terms of tropical optimization}

We start with the general constrained location problem, which is formulated to find all vectors $\bm{x}=(x_{1},x_{2})^{T}\in\mathbb{R}^{2}$ that  
\begin{equation}
\begin{aligned}
&
\text{minimize}
&&
\max_{1\leq j\leq m}(\rho(\bm{x},\bm{r}_{j})+w_{j}),
\\
&
\text{subject to}
&&
\rho(\bm{x},\bm{r}_{j})
\leq
d_{j},
\quad
j=1,\ldots,m;
\\
&&&
s
\leq
x_{1}
\leq
t;
\end{aligned}
\label{P-rhoxrjwj-rhoxrjdj-sx1t}
\end{equation}
where $\bm{r}_{j}=(r_{1j},r_{2j})^{T}\in\mathbb{R}^{2}$ are given vectors and $d_{j},w_{j}\in\mathbb{R}$ with $d_{j}\geq0$ are given numbers for all $j=1,\ldots,m$, and $s,t\in\mathbb{R}$ are given numbers such that $s\leq t$. 

To solve problem \eqref{P-rhoxrjwj-rhoxrjdj-sx1t}, we represent it in terms of the semifield $\mathbb{R}_{\max,+}$ with the maximum in the role of addition and the arithmetic addition in the role of multiplication. Clearly, the context of location analysis guarantees the regularity, in terms of $\mathbb{R}_{\max,+}$, of all vectors involved in the problem setting. 

First, we note that the rectilinear distance between two vectors $\bm{a}=(a_{1},a_{2})^{T}$ and $\bm{b}=(b_{1},b_{2})^{T}$ in terms of the operations in the semifield $\mathbb{R}_{\max,+}$ takes the form
$$
\rho(\bm{a},\bm{b})
=
(a_{1}^{-1}b_{1}\oplus b_{1}^{-1}a_{1})(a_{2}^{-1}b_{2}\oplus b_{2}^{-1}a_{2}).
$$

Consider the distance between the points $\bm{x}$ and $r_{j}$ for each $j=1,\ldots,m$. An application of the above formula and simple algebra give
$$
\rho(\bm{x},\bm{r}_{j})
=
r_{1j}r_{2j}x_{1}^{-1}x_{2}^{-1}
\oplus
r_{1j}^{-1}r_{2j}x_{1}x_{2}^{-1}
\oplus
r_{1j}r_{2j}^{-1}x_{1}^{-1}x_{2}
\oplus
r_{1j}^{-1}r_{2j}^{-1}x_{1}x_{2}.
$$

To represent this distance in a compact vector form, we introduce the vectors $\bm{y}=(y_{1},y_{2})^{T}$ and $\bm{c}_{j}=(c_{1j},c_{2j})^{T}$ for all $j=1,\ldots,m$, with elements
\begin{align*}
y_{1}
&=
x_{1}x_{2},
&
c_{1j}
&=
r_{1j}r_{2j},
\\
y_{2}
&=
x_{1}^{-1}x_{2},
&
c_{2j}
&=
r_{1j}^{-1}r_{2j}.
\end{align*}

It is not difficult to see that the elements of the vector $\bm{x}$ are uniquely determined by those of $\bm{y}$ through the equalities
$$
x_{1}
=
y_{1}^{1/2}y_{2}^{-1/2},
\qquad
x_{2}
=
y_{1}^{1/2}y_{2}^{1/2}.
$$

Note that the above relations between $x_{1},x_{2}$ and $y_{1},y_{2}$ directly correspond to the obvious linear transformations, which are associated in the conventional arithmetic with the interchange between Chebyshev and rectilinear distances.

With the new vector notation, the distance between $\bm{x}$ and $r_{j}$ becomes
$$
\rho(\bm{x},\bm{r}_{j})
=
\bm{y}^{-}\bm{c}_{j}
\oplus
\bm{c}_{j}^{-}\bm{y}.
$$

We are now in a position to rewrite the objective function in problem \eqref{P-rhoxrjwj-rhoxrjdj-sx1t} by using vectors $\bm{p}=(p_{1},p_{2})^{T}$ and $\bm{q}=(q_{1},q_{2})^{T}$ as follows:
$$
\bigoplus_{j=1}^{m}w_{j}(\bm{y}^{-}\bm{c}_{j}\oplus\bm{c}_{j}^{-}\bm{y})
=
\bm{y}^{-}\bm{p}\oplus\bm{q}^{-}\bm{y},
$$
where the right-hand side is obtained by regrouping terms and substitution
$$
\bm{p}
=
\bigoplus_{j=1}^{m}w_{j}\bm{c}_{j},
\qquad
\bm{q}^{-}
=
\bigoplus_{j=1}^{m}w_{j}\bm{c}_{j}^{-}.
$$

Furthermore, we examine the inequality constraints in \eqref{P-rhoxrjwj-rhoxrjdj-sx1t}. The constraints, which involve the distance between vectors, take the form of the inequalities
$$
\bm{y}^{-}\bm{c}_{j}\oplus\bm{c}_{j}^{-}\bm{y}
\leq
d_{j},
\qquad
j=1,\ldots,m.
$$

Note that each inequality is equivalent to the pair of inequalities $\bm{y}^{-}\bm{c}_{j}\leq d_{j}$ and $\bm{c}_{j}^{-}\bm{y}\leq d_{j}$. Consider the first inequality $\bm{y}^{-}\bm{c}_{j}\leq d_{j}$ and verify, using the properties of conjugate transposition, that it is equivalent to the inequality $\bm{c}_{j}\leq d_{j}\bm{y}$. Indeed, multiplication of the former inequality by $\bm{y}$ from the left gives $\bm{c}_{j}\leq\bm{y}\bm{y}^{-}\bm{c}_{j}\leq d_{j}\bm{y}$, whereas the left multiplication of the latter inequality by $\bm{y}^{-}$ results in the former one. The inequality $\bm{c}_{j}^{-}\bm{y}\leq d_{j}$ is equivalent to $\bm{y}\leq d_{j}\bm{c}_{j}$ by similar arguments.

Then, after slight rearrangement of the inequalities obtained, we represent the inequality constraints under consideration in the alternative form
$$
d_{j}^{-1}\bm{c}_{j}
\leq
\bm{y},
\qquad
\bm{y}^{-}
\geq
d_{j}^{-1}\bm{c}_{j}^{-},
\qquad
j=1,\ldots,m.
$$

These inequalities combine to produce the two equivalent inequalities
$$
\bigoplus_{j=1}^{m}d_{j}^{-1}\bm{c}_{j}
\leq
\bm{y},
\qquad
\bm{y}^{-}
\geq
\bigoplus_{j=1}^{m}d_{j}^{-1}\bm{c}_{j}^{-}.
$$

Finally, we replace the last inequalities by the one double inequality
$$
\bm{g}
\leq
\bm{y}
\leq
\bm{h},
$$
where we use the vector notation $\bm{g}=(g_{1},g_{2})^{T}$ and $\bm{h}=(h_{1},h_{2})^{T}$, defined by
$$
\bm{g}
=
\bigoplus_{j=1}^{m}d_{j}^{-1}\bm{c}_{j},
\qquad
\bm{h}^{-}
=
\bigoplus_{j=1}^{m}d_{j}^{-1}\bm{c}_{j}^{-}.
$$

It remains to represent, in terms of the new vector $\bm{y}$, the last inequality constraint
$$
s
\leq
x_{1}
\leq
t.
$$

We rewrite the left and right inequalities as $s^{2}x_{1}^{-1}\leq x_{1}$ and $t^{-2}x_{1}\leq x_{1}^{-1}$, or equivalently, as the inequalities $s^{2}x_{1}^{-1}x_{2}\leq x_{1}x_{2}$ and $t^{-2}x_{1}x_{2}\leq x_{1}^{-1}x_{2}$.

After substitution $y_{1}=x_{1}x_{2}$ and $y_{2}=x_{1}^{-1}x_{2}$, we have the inequalities $s^{2}y_{2}\leq y_{1}$ and $t^{-2}y_{1}\leq y_{2}$. In vector form, these inequalities are given by
$$
\bm{B}\bm{y}
\leq
\bm{y},
$$
where the matrix $\bm{B}$ is defined, using the notation $\mathbb{0}=-\infty$, as follows:
$$
\bm{B}
=
\left(
\begin{array}{cc}
\mathbb{0} & s^{2}
\\
t^{-2} & \mathbb{0}
\end{array}
\right).
$$

Finally, location problem \eqref{P-rhoxrjwj-rhoxrjdj-sx1t} reduces to the tropical optimization problem
\begin{equation}
\begin{aligned}
&
\text{minimize}
&&
\bm{y}^{-}\bm{p}\oplus\bm{q}^{-}\bm{y},
\\
&
\text{subject to}
&&
\bm{B}\bm{y}
\leq
\bm{y},
\\
&&&
\bm{g}
\leq
\bm{y}
\leq
\bm{h},
\end{aligned}
\label{P-ypqy-Byy-gyh}
\end{equation}
which coincides with that of \eqref{P-xpqx-Bxx-gxh}, where the unknown vector $\bm{x}$ is replaced by $\bm{y}$.

\subsection{Derivation of direct solution}

We now apply Theorem~\ref{T-xpqx-Bxx-gxh} to derive a direct solution to problem \eqref{P-ypqy-Byy-gyh}. To describe the results, we need to calculate the matrices
$$
\bm{B}^{\ast}
=
\bm{I}\oplus\bm{B}
=
\left(
\begin{array}{cc}
\mathbb{1} & s^{2}
\\
t^{-2} & \mathbb{1}
\end{array}
\right),
\qquad
\bm{B}^{2}
=
\left(
\begin{array}{cc}
s^{2}t^{-2} & \mathbb{0}
\\
\mathbb{0} & s^{2}t^{-2}
\end{array}
\right)
=
s^{2}t^{-2}\bm{I}.
$$
 
The description also involves a direct representation for the elements of the vectors $\bm{p}=(p_{1},p_{2})^{T}$, $\bm{q}=(q_{1},q_{2})^{T}$, $\bm{g}=(g_{1},g_{2})^{T}$ and $\bm{h}=(h_{1},h_{2})^{T}$, given by
\begin{equation}
p_{i}
=
\bigoplus_{j=1}^{m}w_{j}c_{ij},
\quad
q_{i}^{-}
=
\bigoplus_{j=1}^{m}w_{j}c_{ij}^{-1},
\quad
g_{i}
=
\bigoplus_{j=1}^{m}d_{j}^{-1}c_{ij},
\quad
h_{i}^{-}
=
\bigoplus_{j=1}^{m}d_{j}^{-1}c_{ij}^{-1},
\label{E-pi-qi-gi-hi}
\end{equation}
where $i=1,2$, and $c_{1j}=r_{1j}r_{2j}$ and $c_{2j}=r_{1j}^{-1}r_{2j}$ for all $j=1,\ldots,m$.

We start with calculating, as intermediate results, the row vectors
\begin{align*}
\bm{q}^{-}\bm{B}^{\ast}
&=
\left(
\begin{array}{cc}
q_{1}^{-1}
\oplus
t^{-2}q_{2}^{-1},
&
s^{2}q_{1}^{-1}
\oplus
q_{2}^{-1}
\end{array}
\right),
\\
\bm{h}^{-}\bm{B}^{\ast}
&=
\left(
\begin{array}{cc}
h_{1}^{-1}
\oplus
t^{-2}h_{2}^{-1},
&
s^{2}h_{1}^{-1}
\oplus
h_{2}^{-1}
\end{array}
\right).
\end{align*}

To represent the existence conditions imposed and the minimum value provided by the theorem, we find
\begin{align*}
\bm{q}^{-}\bm{B}^{\ast}\bm{p}
&=
q_{1}^{-1}p_{1}
\oplus
t^{-2}q_{2}^{-1}p_{1}
\oplus
s^{2}q_{1}^{-1}p_{2}
\oplus
q_{2}^{-1}p_{2},
\\
\bm{q}^{-}\bm{B}^{\ast}\bm{g}
&=
q_{1}^{-1}g_{1}
\oplus
t^{-2}q_{2}^{-1}g_{1}
\oplus
s^{2}q_{1}^{-1}g_{2}
\oplus
q_{2}^{-1}g_{2},
\\
\bm{h}^{-}\bm{B}^{\ast}\bm{p}
&=
h_{1}^{-1}p_{1}
\oplus
t^{-2}h_{2}^{-1}p_{1}
\oplus
s^{2}h_{1}^{-1}p_{2}
\oplus
h_{2}^{-1}p_{2},
\\
\bm{h}^{-}\bm{B}^{\ast}\bm{g}
&=
h_{1}^{-1}g_{1}
\oplus
t^{-2}h_{2}^{-1}g_{1}
\oplus
s^{2}h_{1}^{-1}g_{2}
\oplus
h_{2}^{-1}g_{2}.
\end{align*}

According to Theorem~\ref{T-xpqx-Bxx-gxh}, the conditions for problem \eqref{P-ypqy-Byy-gyh} to have solutions are specified by the inequalities
$$
\mathop\mathrm{Tr}(\bm{B})
\leq
\mathbb{1},
\qquad
\bm{h}^{-}\bm{B}^{\ast}\bm{g}
\leq
\mathbb{1}.
$$

Since $\mathop\mathrm{Tr}(\bm{B})=\mathop\mathrm{tr}\bm{B}\oplus\mathop\mathrm{tr}\bm{B}^{2}=s^{2}t^{-2}\leq\mathbb{1}$ if $s\leq t$, the first inequality is obviously valid. The second condition takes the form of the inequality
\begin{equation}
h_{1}^{-1}g_{1}
\oplus
t^{-2}h_{2}^{-1}g_{1}
\oplus
s^{2}h_{1}^{-1}g_{2}
\oplus
h_{2}^{-1}g_{2}
\leq
\mathbb{1}.
\label{I-h1g1-t2h2g1-s2h1g2-h2g2-1}
\end{equation}

An application of \eqref{E-theta-qBp-hBpqBg} to write the minimum value of the objective function yields
\begin{multline}
\theta
=
(q_{1}^{-1}p_{1}
\oplus
t^{-2}q_{2}^{-1}p_{1}
\oplus
s^{2}q_{1}^{-1}p_{2}
\oplus
q_{2}^{-1}p_{2})^{1/2}
\\
\oplus
h_{1}^{-1}p_{1}
\oplus
t^{-2}h_{2}^{-1}p_{1}
\oplus
s^{2}h_{1}^{-1}p_{2}
\oplus
h_{2}^{-1}p_{2}
\\
\oplus
q_{1}^{-1}g_{1}
\oplus
t^{-2}q_{2}^{-1}g_{1}
\oplus
s^{2}q_{1}^{-1}g_{2}
\oplus
q_{2}^{-1}g_{2}.
\label{E-theta-q1p1-t2q2p1-s2q1p2}
\end{multline}

We now describe the solution set of vectors $\bm{y}=(y_{1},y_{2})^{T}$. It follows from Theorem~\ref{T-xpqx-Bxx-gxh} that problem \eqref{P-ypqy-Byy-gyh} has the solution
$$
\bm{y}
=
\bm{B}^{\ast}\bm{u},
$$
where the intermediate vector $\bm{u}=(u_{1},u_{2})^{T}$ satisfies the condition at \eqref{E-gthetap-u-hthetaqB}. 

First, we represent the above vector equality in scalar form as
\begin{equation}
y_{1}
=
u_{1}\oplus s^{2}u_{2},
\qquad
y_{2}
=
t^{-2}u_{1}\oplus u_{2}.
\label{E-y1u1s2u2-y2t2u1u2}
\end{equation}

To describe the set of admissible vectors $\bm{u}$, we take the double inequality at \eqref{E-gthetap-u-hthetaqB} to consider the corresponding scalar inequalities
\begin{align*}
g_{1}\oplus\theta^{-1}p_{1}
&\leq
u_{1}
\leq
(h_{1}^{-1}\oplus\theta^{-1}q_{1}^{-1}\oplus t^{-2}h_{2}^{-1}\oplus\theta^{-1}t^{-2}q_{2}^{-1})^{-1},
\\
g_{2}\oplus\theta^{-1}p_{2}
&\leq
u_{2}
\leq
(s^{2}h_{1}^{-1}\oplus\theta^{-1}s^{2}q_{1}^{-1}\oplus h_{2}^{-1}\oplus\theta^{-1}q_{2}^{-1})^{-1}.
\end{align*}

It is not difficult to verify that at least one of these inequalities holds as an equality. To see this, we can substitute for $\theta$ each term on the right-hand side of \eqref{E-theta-q1p1-t2q2p1-s2q1p2}. Consider, for instance, the case that $\theta=q_{1}^{-1/2}p_{1}^{1/2}$. Under this assumption, we have
\begin{multline*}
p_{1}^{1/2}q_{1}^{1/2}
=
\theta^{-1}p_{1}
\leq
g_{1}\oplus\theta^{-1}p_{1}
\\
\leq
(h_{1}^{-1}\oplus\theta^{-1}q_{1}^{-1}\oplus t^{-2}h_{2}^{-1}\oplus\theta^{-1}t^{-2}q_{2}^{-1})^{-1}
\leq
\theta q_{1}
=
p_{1}^{1/2}q_{1}^{1/2},
\end{multline*}
which means that both left and right parts of the first inequality coincide, and thus this inequality reduces to an equality. The other cases are examined in the same way.

Taking into account that one of the above inequalities acts as an equality, we rewrite them in the one-parametrized form
\begin{equation}
\begin{aligned}
u_{1}
&=
(g_{1}\oplus\theta^{-1}p_{1})^{1-\alpha}
(h_{1}^{-1}\oplus\theta^{-1}q_{1}^{-1}\oplus t^{-2}h_{2}^{-1}\oplus\theta^{-1}t^{-2}q_{2}^{-1})^{-\alpha},
\\
u_{2}
&=
(g_{2}\oplus\theta^{-1}p_{2})^{1-\alpha}
(s^{2}h_{1}^{-1}\oplus\theta^{-1}s^{2}q_{1}^{-1}\oplus h_{2}^{-1}\oplus\theta^{-1}q_{2}^{-1})^{-\alpha},
\end{aligned}
\label{E-u1-u2}
\end{equation}
where $\alpha$ is a real parameter such that $0\leq\alpha\leq1$.

Finally, note that the solution of the initial problem \eqref{P-rhoxrjwj-rhoxrjdj-sx1t} in terms of the vector $\bm{x}=(x_{1},x_{2})^{T}$ can be calculated from the elements of the vector $\bm{y}$ as follows:
\begin{equation}
x_{1}
=
y_{1}^{1/2}y_{2}^{-1/2},
\qquad
x_{2}
=
y_{1}^{1/2}y_{2}^{1/2}.
\label{E-x1y1y2-x2-y1y2}
\end{equation}

\subsection{Direct solutions to location problems}

In this subsection, we turn back to the conventional notation and summarize the result obtained to provide direct, explicit solutions of the general location problem and of its special cases. Consider the general problem formulated in the scalar form
\begin{equation}
\begin{aligned}
&
\text{minimize}
&&
\max_{1\leq j\leq m}(|x_{1}-r_{1j}|+|x_{2}-r_{2j}|+w_{j}),
\\
&
\text{subject to}
&&
|x_{1}-r_{1j}|+|x_{2}-r_{2j}|
\leq
d_{j},
\quad
j=1,\ldots,m;
\\
&&&
s
\leq
x_{1}
\leq
t.
\end{aligned}
\label{P-xqr1jx2r2jwj-xqr1jx2r2dj-sx1t}
\end{equation}

After translating the formulae at \eqref{E-pi-qi-gi-hi}, \eqref{I-h1g1-t2h2g1-s2h1g2-h2g2-1}, \eqref{E-theta-q1p1-t2q2p1-s2q1p2}, \eqref{E-y1u1s2u2-y2t2u1u2}, \eqref{E-u1-u2} and \eqref{E-x1y1y2-x2-y1y2} back into the language of conventional algebra, with eliminating both $y_{1}$ and $y_{2}$, the results of the previous subsection are summarized as follows.
\begin{theorem}
\label{T-xqr1jx2r2jwj-xqr1jx2r2dj-sx1t}
Define the notation
\begin{equation}
\begin{aligned}
p_{1}
&=
\max_{1\leq j\leq m}(w_{j}+r_{1j}+r_{2j}),
&
p_{2}
&=
\max_{1\leq j\leq m}(w_{j}-r_{1j}+r_{2j}),
\\
q_{1}
&=
\min_{1\leq j\leq m}(-w_{j}+r_{1j}+r_{2j}),
&
q_{2}
&=
\min_{1\leq j\leq m}(-w_{j}-r_{1j}+r_{2j}),
\\
g_{1}
&=
\max_{1\leq j\leq m}(-d_{j}+r_{1j}+r_{2j}),
&
g_{2}
&=
\max_{1\leq j\leq m}(-d_{j}-r_{1j}+r_{2j}),
\\
h_{1}
&=
\min_{1\leq j\leq m}(d_{j}+r_{1j}+r_{2j}),
&
h_{2}
&=
\min_{1\leq j\leq m}(d_{j}-r_{1j}+r_{2j}),
\end{aligned}
\label{E-p1p2q1q2g1g2h1h2}
\end{equation}
and suppose that
\begin{equation}
\max(g_{1}-h_{1},
g_{1}-h_{2}-2t,
g_{2}-h_{1}+2s,
g_{2}-h_{2})
\leq
0.
\label{I-g1h1g1h2tg2h1sg2h2}
\end{equation}

Then, the minimum value in problem \eqref{P-xqr1jx2r2jwj-xqr1jx2r2dj-sx1t} is equal to
\begin{multline}
\theta
=
\max(
(p_{1}-q_{1})/2,
(p_{1}-q_{2})/2-t,
(p_{2}-q_{1})/2+s,
(p_{2}-q_{2})/2,
\\
p_{1}-h_{1},
p_{1}-h_{2}-2t,
p_{2}-h_{1}+2s,
p_{2}-h_{2},
\\
g_{1}-q_{1},
g_{1}-q_{2}-2t,
g_{2}-q_{1}+2s,
g_{2}-q_{2}),
\label{E-theta-pqphgq}
\end{multline}
and all solutions are given by
\begin{equation}
\begin{aligned}
x_{1}
&=
\max(u_{1},u_{2}+2s)/2-\max(u_{1}-2t,u_{2})/2,
\\
x_{2}
&=
\max(u_{1},u_{2}+2s)/2+\max(u_{1}-2t,u_{2})/2,
\end{aligned}
\label{E-x1-x2}
\end{equation}
where
\begin{align*}
u_{1}
&=
(1-\alpha)
\max(g_{1},p_{1}-\theta)
+
\alpha
\min(h_{1},q_{1}+\theta,h_{2}+2t,q_{2}+2t+\theta),
\\
u_{2}
&=
(1-\alpha)
\max(g_{2},p_{2}-\theta)
+
\alpha
\min(h_{1}-2s,q_{1}-2s+\theta,h_{2},q_{2}+\theta)
\end{align*}
for all real $\alpha$ such that $0\leq\alpha\leq1$.
\end{theorem}

As a consequence of the theorem, which also takes into account Corollaries~\ref{C-xpqx-Bxx}, \ref{C-xpqx-gxh} and \ref{C-xpqx}, we present solutions to special cases of the problem with reduced sets of constraints and without constraints. Consider the problem, which has only the upper-bound distance constraints. In ordinary notation, the problem is defined as follows:
\begin{equation}
\begin{aligned}
&
\text{minimize}
&&
\max_{1\leq j\leq m}(|x_{1}-r_{1j}|+|x_{2}-r_{2j}|+w_{j}),
\\
&
\text{subject to}
&&
|x_{1}-r_{1j}|+|x_{2}-r_{2j}|
\leq
d_{j},
\quad
j=1,\ldots,m.
\end{aligned}
\label{P-x1r1jx2r2jwj-x1r1jx2r2dj}
\end{equation}

The next statement offers a direct, explicit solution to the problem.
\begin{corollary}
\label{C-x1r1jx2r2jwj-x1r1jx2r2dj}
Under the notation of Theorem~\ref{T-xqr1jx2r2jwj-xqr1jx2r2dj-sx1t}, suppose that
$$
\max(g_{1}-h_{1},g_{2}-h_{2})
\leq
0.
$$

Then, the minimum value in problem \eqref{P-x1r1jx2r2jwj-x1r1jx2r2dj} is equal to
$$
\theta
=
\max((p_{1}-q_{1})/2,(p_{2}-q_{2})/2,p_{1}-h_{1},p_{2}-h_{2},g_{1}-q_{1},g_{2}-q_{2}),
$$
and all solutions are given by
$$
x_{1}
=
(u_{1}-u_{2})/2,
\qquad
x_{2}
=
(u_{1}+u_{2})/2,
$$
where
$$
u_{i}
=
(1-\alpha)\max(g_{i},p_{i}-\theta)
+
\alpha\min(h_{i},q_{i}+\theta),
\qquad
i=1,2,
$$
for all real $\alpha$ such that $0\leq\alpha\leq1$.
\end{corollary}

Furthermore, we examine the problem with the boundary constraints in the form
\begin{equation}
\begin{aligned}
&
\text{minimize}
&&
\max_{1\leq j\leq m}(|x_{1}-r_{1j}|+|x_{2}-r_{2j}|+w_{j}),
\\
&
\text{subject to}
&&
s
\leq
x_{1}
\leq
t.
\end{aligned}
\label{P-x1r1jx2r2jwj-sx1t}
\end{equation}

\begin{corollary}
\label{C-x1r1jx2r2jwj-sx1t}
Under the notation of Theorem~\ref{T-xqr1jx2r2jwj-xqr1jx2r2dj-sx1t}, the minimum value in problem \eqref{P-x1r1jx2r2jwj-sx1t} is equal to
$$
\theta
=
\max(p_{1}-q_{1},p_{1}-q_{2}-2t,p_{2}-q_{1}+2s,p_{2}-q_{2})/2,
$$
and all solutions are given by
\begin{align*}
x_{1}
&=
\max(u_{1},u_{2}+2s)/2-\max(u_{1}-2t,u_{2})/2,
\\
x_{2}
&=
\max(u_{1},u_{2}+2s)/2+\max(u_{1}-2t,u_{2})/2,
\end{align*}
where
\begin{align*}
u_{1}
&=
(1-\alpha)(p_{1}-\theta)
+
\alpha
(\min(q_{1},q_{2}+2t)+\theta),
\\
u_{2}
&=
(1-\alpha)(p_{2}-\theta)
+
\alpha
(\min(q_{1}-2s,q_{2})+\theta)
\end{align*}
for all real $\alpha$ such that $0\leq\alpha\leq1$.
\end{corollary}

Finally, we consider the unconstrained problem
\begin{equation}
\begin{aligned}
&
\text{minimize}
&&
\max_{1\leq j\leq m}(|x_{1}-r_{1j}|+|x_{2}-r_{2j}|+w_{j}).
\end{aligned}
\label{P-x1r1jx2r2jwj}
\end{equation}

A solution to the problem is described as follows.
\begin{corollary}
\label{C-x1r1jx2r2jwj}
Under the notation of Theorem~\ref{T-xqr1jx2r2jwj-xqr1jx2r2dj-sx1t}, the minimum value in problem \eqref{P-x1r1jx2r2jwj} is equal to
$$
\theta
=
\max(p_{1}-q_{1},p_{2}-q_{2})/2,
$$
and all solutions are given by
$$
x_{1}
=
(u_{1}-u_{2})/2,
\qquad
x_{2}
=
(u_{1}+u_{2})/2,
$$
where
$$
u_{i}
=
(1-\alpha)(p_{i}-\theta)
+
\alpha
(q_{i}+\theta),
\qquad
i=1,2,
$$
for all real $\alpha$ such that $0\leq\alpha\leq1$.
\end{corollary}

Note that, after elimination of the intermediate variables $u_{1}$ and $u_{2}$, the solution to the unconstrained problem becomes
\begin{align*}
x_{1}
&=
(1-\alpha)(p_{1}-p_{2})/2
+
\alpha
(q_{1}-q_{2})/2
\\
x_{2}
&=
(2\alpha-1)\theta
+
(1-\alpha)(p_{1}+p_{2})/2
+
\alpha(q_{1}+q_{2})/2,
\end{align*}
which agrees with that derived by geometric \cite{Elzinga1972Geometrical,Francis1972Ageometrical} and algebraic \cite{Krivulin2011Anextremal,Krivulin2015Onanalgebraic} techniques.

\section{Numerical examples and graphical illustrations}
\label{S-NEGI}

We illustrate the results obtained with small artificial examples of optimal location of a facility with respect to $m=3$ given points. The purpose of the illustration is to provide a clear demonstration of the computational technique used, and a transparent graphical representation of the solutions offered. Although the problems under consideration involve only three given points, the examples show strong evidence of the applicability of the method to solve efficiently real-world problems of large scale.

Consider the problem of locating a new point on the plane to minimize the maximum of distances from this point to three given points defined as
$$
\bm{r}_{1}
=
\left(
\begin{array}{c}
1
\\
2
\end{array}
\right),
\qquad
\bm{r}_{2}
=
\left(
\begin{array}{c}
5
\\
9
\end{array}
\right),
\qquad
\bm{r}_{3}
=
\left(
\begin{array}{c}
7
\\
5
\end{array}
\right).
$$

The values of addends corresponding to these points are assumed to be
$$
w_{1}
=
2,
\qquad
w_{2}
=
1,
\qquad
w_{3}
=
1,
$$
whereas the upper bounds on the distances are to be
$$
d_{1}
=
7,
\qquad
d_{2}
=
5,
\qquad
d_{3}
=
5.
$$

Finally, the left and right boundary of the feasible location region are given by 
$$
s
=
4,
\qquad
t
=
8.
$$

To describe the solutions to the problem under various constraints, we first calculate the numbers
\begin{align*}
p_{1}
&=
\max_{1\leq j\leq 3}(w_{j}+r_{1j}+r_{2j})
=
15,
&
p_{2}
&=
\max_{1\leq j\leq 3}(w_{j}-r_{1j}+r_{2j})
=
5,
\\
q_{1}
&=
\min_{1\leq j\leq 3}(-w_{j}+r_{1j}+r_{2j})
=
1,
&
q_{2}
&=
\min_{1\leq j\leq 3}(-w_{j}-r_{1j}+r_{2j})
=
-3,
\\
g_{1}
&=
\max_{1\leq j\leq 3}(-d_{j}+r_{1j}+r_{2j})
=
9,
&
g_{2}
&=
\max_{1\leq j\leq 3}(-d_{j}-r_{1j}+r_{2j})
=
-1,
\\
h_{1}
&=
\min_{1\leq j\leq 3}(d_{j}+r_{1j}+r_{2j})
=
10,
&
h_{2}
&=
\min_{1\leq j\leq 3}(d_{j}-r_{1j}+r_{2j})
=
3.
\end{align*}

We start with problem \eqref{P-x1r1jx2r2jwj} without constraints. According to Corollary~\ref{C-x1r1jx2r2jwj}, we find the minimum
$$
\theta
=
\max(p_{1}-q_{1},p_{2}-q_{2})/2
=
7.
$$

Next, we calculate the intermediate variables
\begin{align*}
u_{1}
&=
(1-\alpha)(p_{1}-\theta)
+
\alpha
(q_{1}+\theta)
=
8,
\\
u_{2}
&=
(1-\alpha)(p_{2}-\theta)
+
\alpha
(q_{2}+\theta)
=
6\alpha-2,
\end{align*}
and finally obtain the solution in the parametrized form
$$
x_{1}
=
(u_{1}-u_{2})/2
=
5-3\alpha,
\qquad
x_{2}
=
(u_{1}+u_{2})/2
=
3+3\alpha,
$$
where $\alpha$ is any real number such that $0\leq\alpha\leq1$.

Substitutions $\alpha=0$ and $\alpha=1$ give two points
$$
\bm{x}^{\prime}
=
\left(
\begin{array}{c}
5
\\
3
\end{array}
\right),
\qquad
\bm{x}^{\prime\prime}
=
\left(
\begin{array}{c}
2
\\
6
\end{array}
\right),
$$
which define the ends of the line segment representing the solutions. 

\begin{figure}[ht]
\setlength{\unitlength}{1mm}
\begin{center}
\begin{picture}(55,58)

\newsavebox\location
\savebox{\location}(55,58)[s]
{
\put(3,6){\circle*{1.5}}

\put(-3,6){\circle{1.5}}
\put(3,12){\circle{1.5}}
\put(9,6){\circle{1.5}}
\put(3,0){\circle{1.5}}

\put(3,6){\line(-1,0){6}}
\put(3,6){\line(0,1){6}}
\put(3,6){\line(1,0){6}}
\put(3,6){\line(0,-1){6}}



\put(15,27){\circle*{1.5}}

\put(12,27){\circle{1.5}}
\put(15,30){\circle{1.5}}
\put(18,27){\circle{1.5}}
\put(15,24){\circle{1.5}}

\put(15,27){\line(-1,0){3}}
\put(15,27){\line(0,1){3}}
\put(15,27){\line(1,0){3}}
\put(15,27){\line(0,-1){3}}



\put(21,15){\circle*{1.5}}

\put(18,15){\circle{1.5}}
\put(21,18){\circle{1.5}}
\put(24,15){\circle{1.5}}
\put(21,12){\circle{1.5}}

\put(21,15){\line(-1,0){3}}
\put(21,15){\line(0,1){3}}
\put(21,15){\line(1,0){3}}
\put(21,15){\line(0,-1){3}}





\put(-6,9){\line(1,1){21}}
\put(-6,9){\line(1,-1){12}}
\put(-6,9){\line(1,0){21}}

\put(27,18){\line(-1,1){12}}
\put(27,18){\line(-1,-1){21}}
\put(27,18){\line(-1,0){21}}

\put(15,9){\circle*{1.5}}
\put(6,18){\circle*{1.5}}
\put(6,18){\thicklines\line(1,-1){9}}
\put(6.1,18.1){\thicklines\line(1,-1){9}}
\put(5.9,17.9){\thicklines\line(1,-1){9}}




\put(5,2){$\bm{r}_{1}$}
\put(18,30){$\bm{r}_{2}$}
\put(24,11){$\bm{r}_{3}$}

\put(18,5){$\bm{x}^{\prime}$}
\put(2,21){$\bm{x}^{\prime\prime}$}

\put(-18,0){\vector(1,0){54}}
\put(0,-15){\vector(0,1){57}}

\put(35,-3){$x_{1}$}
\put(-4,42){$x_{2}$}
}

\put(0,0){\makebox(55,60){\put(18,0){\usebox{\location}}}}

\end{picture}
\hspace{7\unitlength}
\begin{picture}(55,58)

\savebox{\location}(55,58)[s]
{
\put(3,6){\circle*{1.5}}





\put(15,27){\circle*{1.5}}





\put(21,15){\circle*{1.5}}





\put(12,-15){\line(0,1){56}}
\put(24,-15){\line(0,1){56}}


\put(-6,9){\line(1,1){21}}
\put(-6,9){\line(1,-1){12}}

\put(27,18){\line(-1,1){12}}
\put(27,18){\line(-1,-1){21}}


\put(15,9){\circle*{1.5}}
\put(12,12){\circle*{1.5}}
\put(11.7,12.3){\thicklines\line(1,-1){3.6}}
\put(11.8,12.4){\thicklines\line(1,-1){3.6}}
\put(11.6,12.2){\thicklines\line(1,-1){3.6}}



\put(4,3){$\bm{r}_{1}$}
\put(16,24){$\bm{r}_{2}$}
\put(17,16){$\bm{r}_{3}$}

\put(18,5){$\bm{x}^{\prime}$}
\put(7,13){$\bm{x}^{\prime\prime}$}

\put(-18,0){\vector(1,0){54}}
\put(0,-15){\vector(0,1){57}}

\put(35,-3){$x_{1}$}
\put(-4,42){$x_{2}$}

\put(13,-3){$s$}
\put(25,-3){$t$}
}

\put(0,0){\makebox(55,60){\put(18,0){\usebox{\location}}}}
\end{picture}
\end{center}
\caption{Unconstrained solution (left) and solution under boundary constraints (right).}
\label{F-USSUBC}
\end{figure}
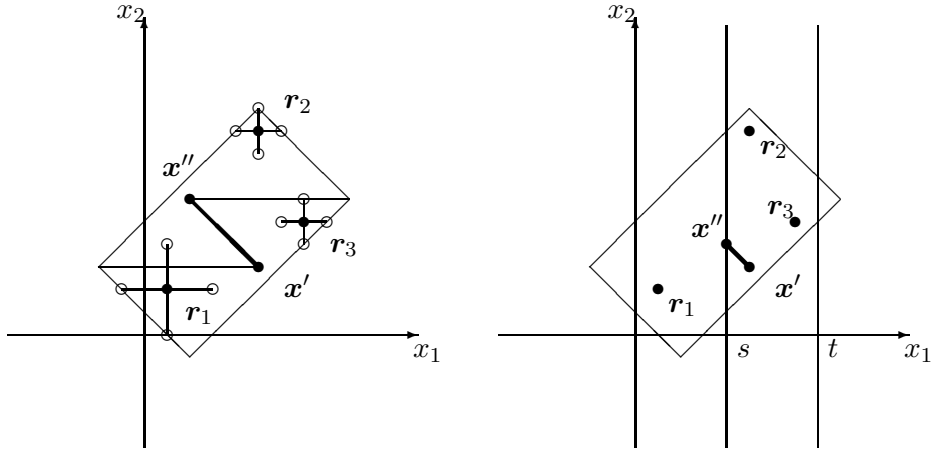

The solution to the unconstrained problem is illustrated in Fig.~\ref{F-USSUBC} (left). The illustration shows the given points $\bm{r}_{1}$, $\bm{r}_{2}$ and $\bm{r}_{3}$, indicated by filled circles. For each point $\bm{r}_{j}$, the four open circles placed distance $w_{j}>0$ from the filled circle designate artificial points to account for the addends. The representation of the solution involves the minimal $45^{\circ}$-tilted rectangle enclosing all artificial points. The solution set is given by the thick line segment, which goes through the centers of the long sides between two horizontal lines drawn through the bottom-left and top-right vertices of the rectangle.

Suppose now that the boundary constraints $s\leq x_{1}\leq t$ are imposed, where $s=4$ and $t=8$. To solve the problem under these constraints, we apply Corollary~\ref{C-x1r1jx2r2jwj-sx1t}.

First, we find that the minimum in the problem
$$
\theta
=
\max(p_{1}-q_{1},p_{1}-q_{2}-2t,p_{2}-q_{1}+2s,p_{2}-q_{2})/2
=
7
$$
remains the same value as for the unconstrained problem examined above.

Furthermore, we calculate the intermediates
\begin{align*}
u_{1}
&=
(1-\alpha)(p_{1}-\theta)
+
\alpha
(\min(q_{1},q_{2}+2t)+\theta)
=
8,
\\
u_{2}
&=
(1-\alpha)(p_{2}-\theta)
+
\alpha
(\min(q_{1}-2s,q_{2})+\theta)
=
-2+2\alpha,
\end{align*}
and then obtain the solution given by
\begin{align*}
x_{1}
&=
\max(u_{1},u_{2}+2s)/2-\max(u_{1}-2t,u_{2})/2
=
5-\alpha,
\\
x_{2}
&=
\max(u_{1},u_{2}+2s)/2+\max(u_{1}-2t,u_{2})/2
=
3+\alpha.
\end{align*}

The solution set forms a line segment with the extreme points, which answer $\alpha=0$ and $\alpha=1$ to be
$$
\bm{x}^{\prime}
=
\left(
\begin{array}{c}
5
\\
3
\end{array}
\right),
\qquad
\bm{x}^{\prime\prime}
=
\left(
\begin{array}{c}
4
\\
4
\end{array}
\right).
$$

The solution is shown in Fig.~\ref{F-USSUBC} (right), where the feasible region is represented as the strip area between two vertical lines at $x_{1}=s$ and $x_{1}=t$.   

Consider the problem, in which the upper-bounds $d_{j}$ on the distances between the new and given points are used instead of boundary constraints examined above. To apply Corollary~\ref{C-x1r1jx2r2jwj-x1r1jx2r2dj}, we verify the condition
$$
\max(g_{1}-h_{1},g_{2}-h_{2})
=
-1
\leq
0.
$$

It follows from the corollary that the minimum in the problem now becomes
$$
\theta
=
\max((p_{1}-q_{1})/2,(p_{2}-q_{2})/2,p_{1}-h_{1},p_{2}-h_{2},g_{1}-q_{1},g_{2}-q_{2})
=
8.
$$

Furthermore, we calculate
\begin{align*}
u_{1}
&=
(1-\alpha)\max(g_{1},p_{1}-\theta)
+
\alpha\min(h_{1},q_{1}+\theta)
=
9,
\\
u_{2}
&=
(1-\alpha)\max(g_{2},p_{2}-\theta)
+
\alpha\min(h_{2},q_{2}+\theta)
=
-1+4\alpha.
\end{align*}

The solution is written as
$$
x_{1}
=
(u_{1}-u_{2})/2
=
5-2\alpha,
\qquad
x_{2}
=
(u_{1}+u_{2})/2
=
4+2\alpha,
$$
and constitutes a line segment having the ends at the points
$$
\bm{x}^{\prime}
=
\left(
\begin{array}{c}
5
\\
4
\end{array}
\right),
\qquad
\bm{x}^{\prime\prime}
=
\left(
\begin{array}{c}
3
\\
6
\end{array}
\right).
$$

\begin{figure}[ht]
\setlength{\unitlength}{1mm}
\begin{center}
\begin{picture}(55,58)

\savebox{\location}(55,58)[s]
{
\put(3,6){\circle*{1.5}}



\put(-18,6){\line(1,1){21}}
\put(-18,6){\line(1,-1){21}}

\put(24,6){\line(-1,1){21}}
\put(24,6){\line(-1,-1){21}}

\put(15,27){\circle*{1.5}}



\put(0,27){\line(1,1){15}}
\put(0,27){\line(1,-1){15}}

\put(30,27){\line(-1,1){15}}
\put(30,27){\line(-1,-1){15}}

\put(21,15){\circle*{1.5}}



\put(6,15){\line(1,1){15}}
\put(6,15){\line(1,-1){15}}

\put(36,15){\line(-1,1){15}}
\put(36,15){\line(-1,-1){15}}


\put(-6,9){\line(1,1){21}}
\put(-6,9){\line(1,-1){12}}

\put(27,18){\line(-1,1){12}}
\put(27,18){\line(-1,-1){21}}



\put(15,12){\circle*{1.5}}
\put(9,18){\circle*{1.5}}
\put(9,18){\thicklines\line(1,-1){6}}
\put(9.1,18.1){\thicklines\line(1,-1){6}}
\put(8.9,17.9){\thicklines\line(1,-1){6}}


\put(4,3){$\bm{r}_{1}$}
\put(10,29){$\bm{r}_{2}$}
\put(24,12){$\bm{r}_{3}$}

\put(15,9){$\bm{x}^{\prime}$}
\put(4,17){$\bm{x}^{\prime\prime}$}

\put(-18,0){\vector(1,0){54}}
\put(0,-15){\vector(0,1){57}}

\put(35,-3){$x_{1}$}
\put(-4,42){$x_{2}$}
}

\put(0,0){\makebox(55,60){\put(18,0){\usebox{\location}}}}

\end{picture}
\hspace{7\unitlength}
\begin{picture}(55,58)

\savebox{\location}(55,58)[s]
{
\put(3,6){\circle*{1.5}}



\put(-18,6){\line(1,1){21}}
\put(-18,6){\line(1,-1){21}}

\put(24,6){\line(-1,1){21}}
\put(24,6){\line(-1,-1){21}}

\put(15,27){\circle*{1.5}}



\put(0,27){\line(1,1){15}}
\put(0,27){\line(1,-1){15}}

\put(30,27){\line(-1,1){15}}
\put(30,27){\line(-1,-1){15}}

\put(21,15){\circle*{1.5}}



\put(6,15){\line(1,1){15}}
\put(6,15){\line(1,-1){15}}

\put(36,15){\line(-1,1){15}}
\put(36,15){\line(-1,-1){15}}

\put(12,-15){\line(0,1){56}}
\put(24,-15){\line(0,1){56}}

\put(-6,9){\line(1,1){21}}
\put(-6,9){\line(1,-1){12}}

\put(27,18){\line(-1,1){12}}
\put(27,18){\line(-1,-1){21}}




\put(15,12){\circle*{1.5}}
\put(12,15){\circle*{1.5}}
\put(11.7,15.3){\thicklines\line(1,-1){3.6}}
\put(11.8,15.4){\thicklines\line(1,-1){3.6}}
\put(11.6,15.2){\thicklines\line(1,-1){3.6}}

\put(4,3){$\bm{r}_{1}$}
\put(17,30){$\bm{r}_{2}$}
\put(25,12){$\bm{r}_{3}$}

\put(15,9){$\bm{x}^{\prime}$}
\put(14,18){$\bm{x}^{\prime\prime}$}

\put(-18,0){\vector(1,0){54}}
\put(0,-15){\vector(0,1){57}}

\put(35,-3){$x_{1}$}
\put(-4,42){$x_{2}$}

\put(13,-3){$s$}
\put(25,-3){$t$}
}

\put(0,0){\makebox(55,60){\put(18,0){\usebox{\location}}}}

\end{picture}
\end{center}
\caption{Solutions under upper-bound distance constraints (left), and under both boundary and upper-bound distance constraints (right).}
\label{F-SUUBDCUBBUBC}
\end{figure}
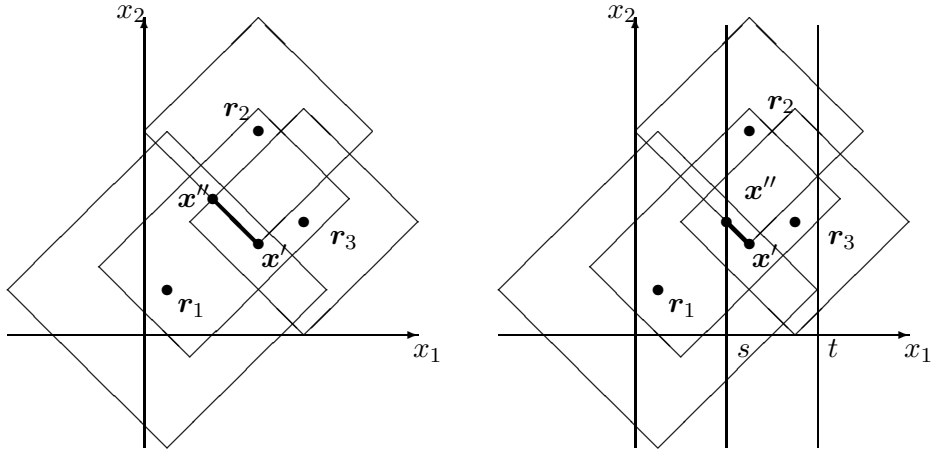

A graphical illustration of the solution is provided in Fig.~\ref{F-SUUBDCUBBUBC} (left). The plot demonstrates the feasible location area as the intersection of turned squares drawn for each point $\bm{r}_{j}$ to have the distance from the vertices of the square to the point equal to $d_{j}$. A thick line segment that coincides with the lower long side of the small turned rectangle, which represents the feasible area, shows the solution.

We conclude this section with the solution to the general location problem, which combines both boundary and upper-bound distance constraints. By following the solution offered by Theorem~\ref{T-xqr1jx2r2jwj-xqr1jx2r2dj-sx1t}, we begin with the validation of the condition
$$
\max(g_{1}-h_{1},
g_{1}-h_{2}-2t,
g_{2}-h_{1}+2s,
g_{2}-h_{2})
=
-10
\leq
0.
$$

The evaluation of the minimum value yields
\begin{multline*}
\theta
=
\max(
(p_{1}-q_{1})/2,
(p_{1}-q_{2})/2-t,
(p_{2}-q_{1})/2+s,
(p_{2}-q_{2})/2,
\\
p_{1}-h_{1},
p_{1}-h_{2}-2t,
p_{2}-h_{1}+2s,
p_{2}-h_{2},
\\
g_{1}-q_{1},
g_{1}-q_{2}-2t,
g_{2}-q_{1}+2s,
g_{2}-q_{2})
=
8.
\end{multline*}

After calculation of the intermediate expressions
\begin{align*}
u_{1}
&=
(1-\alpha)
\max(g_{1},p_{1}-\theta)
+
\alpha
\min(h_{1},q_{1}+\theta,h_{2}+2t,q_{2}+2t+\theta)
=
9,
\\
u_{2}
&=
(1-\alpha)
\max(g_{2},p_{2}-\theta)
+
\alpha
\min(h_{1}-2s,q_{1}-2s+\theta,h_{2},q_{2}+\theta)
=
-1+2\alpha,
\end{align*}
we obtain the solution given by
\begin{align*}
x_{1}
&=
\max(u_{1},u_{2}+2s)/2-\max(u_{1}-2t,u_{2})/2,
=
5-\alpha,
\\
x_{2}
&=
\max(u_{1},u_{2}+2s)/2+\max(u_{1}-2t,u_{2})/2
=
4+\alpha.
\end{align*}

The solution to the problem is depicted in Fig.~\ref{F-SUUBDCUBBUBC} (right) by the thick line segment between the points
$$
\bm{x}^{\prime}
=
\left(
\begin{array}{c}
5
\\
4
\end{array}
\right),
\qquad
\bm{x}^{\prime\prime}
=
\left(
\begin{array}{c}
4
\\
5
\end{array}
\right),
$$
which correspond to setting $\alpha=0$ and $\alpha=1$.

\section{Application to CCTV monitoring facility location}
\label{S-ACCTVMFL}

In this section, we present an application to a real-world problem that arises in the deployment of CCTV video surveillance systems in the indoor environment, including office, industrial, commercial, educational, social, health-care and other buildings. 

A typical CCTV system is composed of three major components \cite{Cieszynski2004Closed,Cohen2009CCTV,SPAWAR2013CCTV}: imaging sensors (video cameras) generating an input video stream, a transmission system to transmit video signal data, and a central video monitoring/processing facility. The design and deployment of CCTV systems in the indoor environment give rise to a range of location problems. Specifically, the problems of camera placement consist in finding optimal locations of cameras in a surveillance zone under various operational objectives and constraints, such as to maximize the coverage area subject to a fixed number of available cameras (see, e.g., overviews in \cite{Zhao2013Approximate,Liu2016Recent}). Below, we examine a different problem with the assumption that the placement of cameras in a CCTV system is already fixed and the problem is to determine the optimal location of the central monitoring facility to reduce losses in the wired transmission network of the system under some technological constraints. 

\begin{figure}[ht]
\setlength{\unitlength}{1mm}
\begin{center}
\begin{picture}(120,145)

\newsavebox\emptyfloorx
\savebox{\emptyfloorx}(72,36)[s]
{\put(0,0){\makebox(72,36)[s]
{\pmb{\pmb{\thicklines
\put(0,20){\line(2,1){32}}
\put(32,36){\line(2,-1){40}}
\put(0,20){\line(2,-1){40}}
\put(40,0){\line(2,1){32}}
}}}}
}

\newsavebox\emptyfloory
\savebox{\emptyfloory}(120,60)[s]
{\put(0,0){\makebox(120,60)[s]
{\pmb{\pmb{\thicklines
\put(0,20){\line(2,1){48}}
\put(80,60){\line(2,-1){40}}
\put(0,20){\line(2,-1){40}}
\put(40,0){\line(2,1){80}}
}}}}
}

\newsavebox\emptyhalffloorx
\savebox{\emptyhalffloorx}(72,36)[s]
{\put(0,0){\makebox(72,36)[s]
{\pmb{\pmb{\thicklines
\put(0,20){\line(2,-1){24}}
\put(0,20){\line(0,1){15}}
}}}}
}

\newsavebox\emptyhalffloory
\savebox{\emptyhalffloory}(120,55)[s]
{\put(0,0){\makebox(120,60)[s]
{\pmb{\pmb{\thicklines
\put(0,20){\line(2,-1){40}}
\put(40,0){\line(2,1){80}}
\put(0,20){\line(0,1){15}}
\put(40,0){\line(0,1){15}}
\put(120,40){\line(0,1){15}}
}}}}
}

\newsavebox\floorx
\savebox{\floorx}(72,36)[s]
{\put(0,0){\makebox(72,36)[s]{\usebox{\emptyfloorx}}}
\pmb{\thicklines
\multiput(0,20)(13.4,-6.7){3}{\line(2,1){24}}
\multiput(24,32)(13.4,-6.7){3}{\line(2,-1){8}}
\multiput(24,32)(13.4,-6.7){3}{\line(0,-1){13}}
\multiput(32,28)(13.4,-6.7){3}{\line(0,-1){5}}
}
}

\newsavebox\floory
\savebox{\floory}(120,60)[s]
{\put(0,0){\makebox(120,60)[s]{\usebox{\emptyfloory}}}
\pmb{\thicklines
\multiput(0,20)(16,8){5}{\line(2,-1){10}}
\multiput(16,28)(16,8){5}{\line(0,-1){10}}
\multiput(10,15)(16,8){5}{\line(2,1){10}}
\multiput(26,23)(16,8){4}{\line(0,-1){8.5}}
\multiput(20,20)(16,8){5}{\line(0,-1){10}}
%
\multiput(18.5,10.5)(16,8){5}{\line(2,-1){21.5}}
\multiput(24,13.5)(16,8){5}{\line(0,-1){5.5}}
\multiput(34,18.5)(16,8){5}{\line(0,-1){15.5}}
\multiput(24,13.5)(16,8){5}{\line(2,1){10}}
}
}

\newsavebox\flooryzero
\savebox{\flooryzero}(120,64)[s]
{\put(0,0){\makebox(120,60)[s]{\usebox{\emptyfloory}}}
{\pmb{\pmb{\thicklines
\put(48,44){\line(2,1){32}}
\put(48,44){\line(-2,1){4}}
\put(80,60){\line(-2,1){4}}
}}}
\multiput(0,20)(1,0.5){80}{\line(1,0){2}}
\multiput(0,20)(1,-0.5){40}{\line(1,0){2}}
\multiput(120,40)(-1,0.5){40}{\line(-1,0){2}}
\multiput(120,40)(-1,-0.5){80}{\line(-1,0){2}}
}

\newsavebox\floorxone
\savebox{\floorxone}(72,36)[s]
{\put(0,0){\makebox(72,36)[s]{\usebox{\floorx}}}
\multiput(5,19)(13.4,-6.7){3}{\faVideoCamera}
}

\newsavebox\flooryone
\savebox{\flooryone}(120,60)[s]
{\put(0,0){\makebox(120,60)[s]{\usebox{\floory}}}
\multiput(5,19)(16,8){2}{\faVideoCamera}
\multiput(53,43)(16,8){2}{\faVideoCamera}
\multiput(46,9)(16,8){2}{\rotatebox{180}{\faVideoCamera}}
\multiput(112,41)(16,8){1}{\rotatebox{180}{\faVideoCamera}}
}

\put(8,109){\usebox{\floorxone}}
\put(8,93){\usebox{\emptyhalffloorx}}

\put(0,65){\usebox{\flooryone}}
\put(0,50){\usebox{\emptyhalffloory}}

\put(0,5){\usebox{\flooryzero}}

{
\pmb{\pmb{\pmb{\thicklines\multiput(77,42)(0,3){7}{\line(0,1){2}}}}}
}

\put(71,45){{\Large\faLaptop}\ \ {\large\faDatabase}}

\put(5,28){\line(0,1){57}}
\put(5,28){\line(2,-1){20}}
\put(25,18){\line(2,1){50}}

\put(114,49){\line(0,1){56}}
\put(114,49){\line(-2,-1){24}}
\put(90,37){\line(-2,1){12}}

\put(93,122){\rotatebox{-27}{$N$th Floor}}

\put(93,63){\rotatebox{-27}{Ground Floor}}

\put(0,26){\vector(2,-1){50}}
\put(52,1){$x_{1}$}

\put(0,26){\vector(2,1){90}}
\put(92,70){$x_{2}$}

\put(12,27){$(r_{1j},r_{2j})$}
\put(8,50){$w_{j}$}

\put(0,22){$s$}
\put(40,1){$t$}

\end{picture}
\end{center}
\caption{Location of CCTV monitoring facility in multi-floor building environment.}
\label{F-CCTV}
\end{figure}
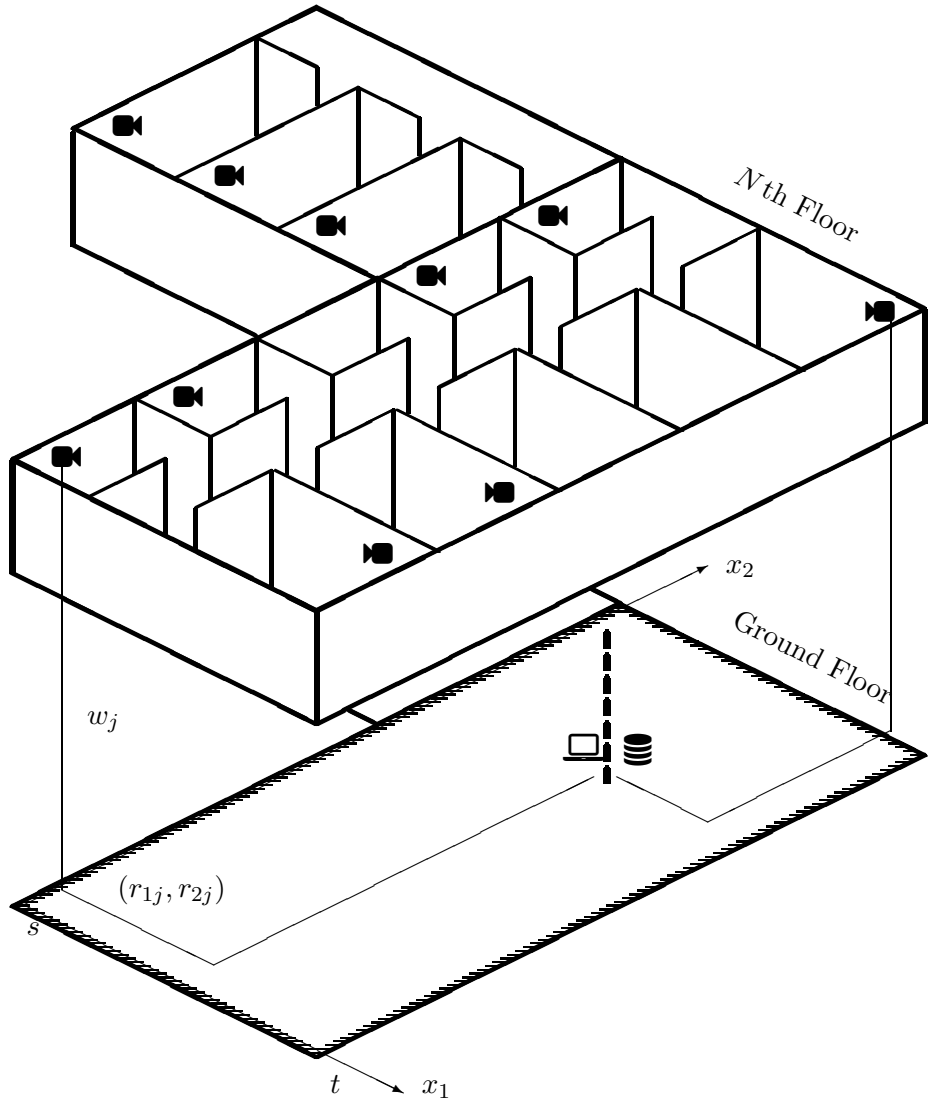

We consider the CCTV video surveillance system, which is set up in a multi-floor building that is composed by rectangular shapes, with rectangular rooms and corridors at each floor, as illustrated by the scheme in Fig.~\ref{F-CCTV}. The intra-building conduit system consists of vertical riser shafts, horizontal cable trays, ladder racks and other facilities to provide full connectivity between any points in the building through the paths, which are parallel or perpendicular to the walls and to the ceiling. Rooms on all floors are equipped with surveillance cameras mounted at upper corners, where two walls and the ceiling meet at a right angle. To transmit video data, every camera is directly connected by a coaxial cable using intra-building cable runs to a central control viewing room, located in a dedicated area on the ground floor, which accommodates monitoring, data storage and video analytics facilities. In the scheme in Fig.~\ref{F-CCTV}, the feasible location region is surrounded by hatched border.  

The use of coaxial wires imposes constraints on the maximum cable distance, which ranges, depending on the gauge of the cable, from several tens to a few hundreds of metres \cite{Cohen2009CCTV,SPAWAR2013CCTV}. Another critical issue is the increasing attenuation (loss) of video signal as the transmission distance and frequency increase, which makes it difficult to handle video data, especially when transmitting high-quality video. Since the length of the wire between each camera and the control room depends on where the room is located, the attenuation is reduced by appropriate room location. Considering that the signal attenuation is measured proportional to the transmission distance, the attenuation can be identified with the distance to formulate the following minimax problem. 

For a given placement of $m$ surveillance cameras in a multi-floor building environment, we need to find the optimal location of the central monitoring facility (the control room) in a feasible area on the ground floor to minimize the maximum attenuation in the wired transmission network subject to maximum distance constraints on the wires between cameras and the facility. As it is easy to see, the problem takes the form of \eqref{P-xqr1jx2r2jwj-xqr1jx2r2dj-sx1t} and, therefore, has the direct solution provided by Theorem~\ref{T-xqr1jx2r2jwj-xqr1jx2r2dj-sx1t}.

The procedure offered by the theorem involves simple straightforward calculations using the horizontal coordinates $(r_{1j},r_{2j})$ and the vertical height $w_{j}$ of all cameras $j=1,\ldots,m$. Given a common distance constraint $d$ on the wire length, the individual constraint for camera $j$ is calculated as $d_{j}=d-w_{j}$. The most computationally intensive part of the procedure is the evaluation of the intermediate variables $p_{i}$, $q_{i}$, $g_{i}$ and $h_{i}$ for $i=1,2$ according to \eqref{E-p1p2q1q2g1g2h1h2}, which takes at most linear time with respect to the number of cameras. The condition \eqref{I-g1h1g1h2tg2h1sg2h2} is verified to ensure that the distance bounds $d_{j}$ and the location bounds $s$ and $t$ can be simultaneously satisfied to provide non-empty feasible location region. The procedure completes with the evaluation of the optimal distance $\theta$ using \eqref{E-theta-pqphgq} and the derivation of parametrized representation for the coordinates $(x_{1},x_{2})$ of the optimal location zone in the form of \eqref{E-x1-x2}.

Solutions to the central monitoring facility location problem are then obtained by plotting the optimal location area, which typically takes the form of a line segment (depicted in Fig.~\ref{F-CCTV} by a thick dashed line), on the ground floor map to determine an appropriate place to deploy the facility and to develop the transmission network.

\section{Conclusions}
\label{S-C}
In this paper, we used tropical mathematics to derive new solutions to constrained minimax single-facility location problems with addends on the plane with rectilinear distance. Tropical mathematics, which deals with the theory and application of algebraic systems with idempotent addition, offers a useful problem formulation and solution framework to provide direct, explicit solutions to a range of classical and novel problems in operations research, management science and other fields. 

To handle the location problems under study, we have formulated these problems in the tropical mathematics setting, and then solved them by applying recent results in tropical optimization. In contrast to the known solution approaches that are mainly based on numerical iterative algorithms, including linear programming and graph optimization, the solutions obtained are given in a simple closed form, which is ready for both further analysis by analytical techniques and straightforward computation with no more than a linear computational cost. As the solution method for the location problems, we have proposed an explicit computational technique that involves easy direct computations, and may help to augment and supersede known approaches when indirect algorithmic solutions are less preferred or even impossible.

Future research will focus on the solution of the problem with additional constraints to define a more general feasible location area on the plane. Extensions of the approach to solve other minimax location problems, including rectilinear problems in three-dimensional space and multi-facility location problems, are also of interest.

\section*{Acknowledgements}
This work was supported in part by the Russian Foundation for Humanities (grant No. 16-02-00059). The author thanks three referees for valuable comments and suggestions, which have been incorporated into the final version of the manuscript.

\bibliographystyle{utphys}

\bibliography{Using_tropical_optimization_to_solve_constrained_minimax_single-facility_location_problems_with_rectilinear_distance}

\end{document}